\documentclass[reqno, 12pts]{amsart}
\usepackage{amsmath,amssymb,amsthm,amsfonts,graphicx,color}
\usepackage[colorlinks=true]{hyperref}
\hypersetup{urlcolor=blue, citecolor=blue}

\usepackage{textcomp}
\usepackage{yfonts}
\usepackage{amssymb}
\usepackage[hmargin=3cm, vmargin=3cm]{geometry}
\usepackage{hyperref}

\def \D{\Delta}
\def\R{\mathbb R}
\def\N{\mathbb N}

\def\ge{\geq}
\def\De{\Delta} 

\def\i0i{\int_0^\infty}

\title[Liouville type theorems for solutions of the weighted fractional Lane-Emden system ]{Liouville type theorems for stable solutions of the weighted fractional Lane-Emden system}
\author{ Hatem Hajlaoui }
\address{ (Hatem Hajlaoui) Institut Sup\'erieur des Math\'ematiques Appliqu\'ees et de l'Informatique. Universit\'e de Kairouan, Tunisie.
\newline    }
\email{hajlaouihatem@gmail.com}

\numberwithin{equation}{section}
\newtheorem{theorem}{Theorem}[section]
\newtheorem{lemma}[theorem]{Lemma}
\newtheorem{definition}[theorem]{Definition}
\newtheorem{rem}{Remark}
\newtheorem{cor}[theorem]{Corollary}
\newtheorem{proposition}{Proposition}[section]
\allowdisplaybreaks

\let\ssection=\section\renewcommand{\section}{\setcounter{equation}{0}\ssection}

\begin{document}
\maketitle

\begin{abstract}In this paper, we prove Liouville type theorems for  stable  solutions to the weighted fractional Lane-Emden system
		\begin{align*}
 (-\Delta)^s u = h(x)v^p,\quad (-\Delta)^s v= h(x)u^q, \quad u,v>0\quad \mbox{in }\;\mathbb{R}^N,
\end{align*}
  where $1<q\leq p$ and
 $h$ is a positive continuous function in $\R^N$ satisfying
 $\displaystyle{\liminf_{|x|\to \infty}}\frac{h(x)}{|x|^\ell} > 0$ with $\ell > 0.$ Our results generalize the results established in \cite{HHM16} for the Laplacian case (correspond to $s=1$) and improve the previous work \cite{TuanHoang21}. As a consequence, we prove classification result for stable solutions to the weighted fractional Lane-Emden equation
$(-\Delta)^s u = h(x)u^p$ in $\mathbb{R}^N$.

\end{abstract}

\noindent
{\it \footnotesize 2010 Mathematics Subject Classification:} {\scriptsize 35J55, 35J65, 35B33, 35B65.}\\
{\it \footnotesize Key words: Liouville type theorems, Stable solutions, Weighted fractional Lane-Emden system and equation} {\scriptsize }

\tableofcontents

\section{Introduction and main results}

Let $s\in (0,\,1)$ and consider the weighted fractional Lane-Emden system
\begin{align}\label{1.1}
 (-\Delta)^s u = h(x)v^p,\quad (-\Delta)^s v= h(x)u^q, \quad u,v>0\quad \mbox{in }\;\mathbb{R}^N,
\end{align}
  where $p, q \geq 1$ and $h$ is a function satisfying
  $$ h\in C(\mathbb{R}^N),\quad h>0, \quad \liminf_{|x|\to \infty}\frac{h(x)}{|x|^\ell} > 0,  \; \ell > 0.
$$
Assume that $u \in C^{2\alpha}_{loc} (\mathbb R^N)\cap \mathcal L_s(\R^N)$ for some $\alpha >s$ with
\[
\mathcal L_s(\R^N)= \left\{u :\R^N \to \R; \, \int_{\R^N} \frac{|u(y)|}{(1 + |y|)^{N+2s}} dy < \infty \right\}.
\]
Then the fractional Laplacian is defined by
\begin{align*}
(-\De)^s u(x) = c_{N,s} P.V.\int_{\R^N} \frac{u(x) - u(y)}{|x-y|^{N+2s}} dy,
\end{align*}
where $c_{N,s}$ is a normalization constant. Clearly,  $(-\De)^s u(x)$ is well-defined at any $x\in\mathbb R^N,$ for any $u$ in $C^{2\alpha}_{loc} (\mathbb R^N)\cap \mathcal L_s(\R^N)$ for some $\alpha >s.$ In what follows, all solutions of \eqref{1.1} are considered in the above space.

\medskip
In this paper, we are interested in the stable solutions to system
\eqref{1.1}. Motivated by \cite{Mon05,Cow13,FG13}, we adopt the following definition of stability
\begin{definition}\label{d.1.1}
A solution $(u,v)$ of \eqref{1.1} is said to be stable if there are two positive functions $\zeta_1, \zeta_2 \in C^{2\alpha}_{loc} (\mathbb R^N)\cap \mathcal L_s(\R^N)$ ($\alpha >s$) satisfying
\begin{align}\label{1.2}
 (-\Delta)^s \zeta_1 = phv^{p-1}\zeta_2,\quad (-\Delta)^s \zeta_2 = qhu^{q-1}\zeta_1.
\end{align}
\end{definition}

\medskip
With $s = 1$ in \eqref{1.1}, the problem of nonexistence of stable solution has been extensively studied by many authors in the past years, see for instance \cite{ WY13, Cow13, CDG14, Hu15, HZ16, Hu18, HHM16, MY19}.
If moreover $h\equiv 1$, the system is reduced to the classical Lane Emden system
\begin{align}\label{1.3}
 -\Delta u = v^p,\quad -\Delta v= u^q, \quad u,v>0\quad \mbox{in }\;\mathbb{R}^N.
\end{align}
It has been conjectured that \eqref{1.3} admits smooth solutions if and only if $p, q>0$ and
$$\frac{1}{p+1}+\frac{1}{q+1} \leq 1-\frac{2}{N}.$$
This conjecture was completely solved for the radial case by \cite{Mit96} (see also \cite{SZ96}) and for general case in low dimensions $N\leq 4$, see \cite{Sou09}. Very recently, Mtiri-Ye \cite{MY19} proved this conjecture for classical solutions which are stable at infinity.

\medskip
For $p, q \geq 1$, Chen-Dupaigne-Ghergu established  in \cite{CDG14} the optimal Liouville type result for radial stable solutions of \eqref{1.3}. In \cite{Cow13}, Cowan developed a new approach to deal with general stable solutions (radial or not). Let
\begin{align}\label{1.4}
t_{0}^{\pm} = \sqrt{\frac{pq(q+1)}{p+1}}\pm\sqrt{\frac{pq(q+1)}{p+1}-\sqrt{\frac{pq(q+1)}{p+1}}},
\end{align}
by combining integral estimates derived from the stability, comparison principle and a bootstrap argument, he proved in \cite{Cow13} the nonexistence result for stable solutions of \eqref{1.3} for $p \geq q > 2t_0^-.$

\medskip
The new approach of Cowan was exploited by many authors to various elliptic systems \cite{Hu15,Hu18,HZ16,HHM16,DP17,Duong17}. In \cite{Hu15}, Hu adapted this approach to study stable solutions of \begin{align}\label{1.5}
 -\Delta u = h(x)v^p,\quad -\Delta v= h(x)u^q, \quad u,v>0\quad \mbox{in }\;\mathbb{R}^N,
\end{align} with $h(x)=(1+|x|^2)^\frac{\ell}{2}$, $\ell\geq 0$, obtaining then a generalization of the results in \cite{Cow13}.  Later on, the author of the present paper with his collaborators \cite{HHM16} improved the results in \cite{Hu15}. In particular, they established a new inverse comparison principle (analogue to \eqref{2.8} below) which is crucial to handle the cases $1 < q \leq \frac{4}{3}.$ More precisely, they obtained the same classification result provided by Theorem \ref{main1} (see below) of stable solutions to \eqref{1.5} with radial function $h$ satisfying $h(x)\geq C (1+|x|^2)^\frac{\ell}{2}$, $\ell\geq 0$.

\medskip
For the non local Lane-Emden system \eqref{1.1}, the first result (up to our knowledge) classifying stable solutions is proved in \cite{TuanHoang21} for $h = 1$ and $s \in (0, 1)$:
\begin{align}\label{1.6}
 (-\Delta)^s u = v^p,\quad (-\Delta)^s v= u^q, \quad u,v>0\quad \mbox{in }\;\mathbb{R}^N.
\end{align}
Adopting the approach in \cite{Cow13}, Duong-Nguyen established the following nonexistence result.
\begingroup
\setcounter{theorem}{0} 
\renewcommand\thetheorem{\Alph{theorem}}
\begin{theorem}
\label{Duong} Let $0<s<1$.
\begin{enumerate}
	\item 	If  $p\geq q>\frac{4}{3}$ and
	\begin{equation}\label{1.7}
	N<2s+\frac{4s(p+1)}{pq-1}t_0^+,
	\end{equation}
	then the system \eqref{1.6} has no stable solution.
	\item  If $1<q\leq \min (p,\frac{4}{3})$, $t_0^-<\frac{q}{2}$ and \eqref{1.7}, then \eqref{1.6} has no stable solution.
\end{enumerate}
\end{theorem}
\endgroup

In this paper, our aim is to generalize \cite{HHM16, TuanHoang21}. Our main results state as follows
\begin{theorem}\label{main1}Suppose that $h$ satisfies $(H)$ and let $x_0$ be the largest root of the polynomial $H(x) = $
\begin{align}\label{H}
 x^4 -
\frac{16pq(q+1)(p+1)}{(pq-1)^2}x^2 +
\frac{16pq(q+1)(p+1)(p+q+2)}{(pq-1)^3}x
-\frac{16pq(q+1)^2(p+1)^2}{(pq-1)^4}.
\end{align}
 \begin{enumerate}
\item
 If $\frac{4}{3}< q \leq p$ then \eqref{1.1} has no stable solution if $N<2s+(2s+\ell)x_0.$
 In particular, if $N \leq 10s+4\ell,$ \eqref{1.1} has no stable
solution for all $\frac{4}{3}<  q \leq p$.
\item If $1<q\leq \min(\frac{4}{3}, p)$, then \eqref{1.1} has no bounded stable solution, if
\begin{equation}\label{1.9}
  N < 2s + \left[\frac{q}{2}+\frac{(2-q)(p q -1)}{(p+q-2)(p+1)}\right](2s+\ell)x_0.
\end{equation}
\end{enumerate}
Therefore, if $N \leq 6s+2\ell,$ the system \eqref{1.1} has no bounded stable
solution for all $p\geq q > 1.$
\end{theorem}
\begin{itemize}
  \item By \cite[Lemma 6]{HHM16},  for any $1<q\leq p$, we have $2t_0^+\frac{p+1}{pq-1}\leq x_0$, and the equality holds if and only if $p=q$. Hence the range of nonexistence result in Theorem \ref{main1} is larger than that provided by Theorem \ref{Duong}.
  \item Using \cite[Remark 2]{HHM16}, we have $2t_0^-<q$ if $q> \frac43.$ It means that the classification result in Theorem \ref{Duong} is valid for $\frac43 <q\leq p.$ However, our approach allows us to provide, for the first time, a rigorous proof for the nonexistence of stable solution to the system \eqref{1.1} with $1<q\leq \frac43.$
 \end{itemize}

\medskip
Consider now the weighted fractional Lane Emden equation:
\begin{align}\label{1.10}
 (-\Delta)^s u = h(x)|u|^{p-1}u\quad \mbox{in }\;\mathbb{R}^N.
\end{align}
For the local case ($s=1$) and when $h\equiv 1,$ Farina completely classified in \cite{Far07} finite Morse index solutions to
\begin{align}\label{1.11}
 -\Delta u = |u|^{p-1}u\quad \mbox{in }\;\mathbb{R}^N.
\end{align}
 He proved that \eqref{1.11} has nontrivial classical solution with finite Morse index if and only if $N\geq3$, $p=\frac{N+2}{N-2}$ or $N\geq11$ and $p\geq p_{JL}=p_{JL}(N,0),$ where $p_{JL}$ stands for the Joseph-Lundgren exponent \cite{jl} (see also \cite{gnw}). Later on, Dancer-Du-Guo \cite{DDG} obtained a sharp critical exponent $p_{JL}(N,\ell)$ with respect to the existence of nontrivial stable solutions $u\in W_{loc}^{1,2}(\mathbb{R}^N)\cap L_{loc}^{\infty}(\mathbb{R}^N)$ of \eqref{1.10} with $s=1$, $h(x)=|x|^\ell,$ $\ell>-2.$ Here $p_{JL}(N,\ell)$ is given by
$$
p_{JL}(N,\ell)=
\left\{
  \begin{array}{ll}
    +\infty & \hbox{if}\quad 1\leq N\leq 10+4\ell, \\
    \frac{(N-2)^2-2(\ell+2)(\ell+N)-2\sqrt{(\ell+2)^3(\ell+2N-2)}}{(N-4\ell-10)(N-2)} & \hbox{if}\quad  N> 10+4\ell.
  \end{array}
\right.
$$
Finally, the same resultat is proved to be true without the local boundedness assumption by Wang-Ye \cite{WY12}.
Many other papers studied stable solutions of \eqref{1.10} with $s=1$, see for instance \cite{CF12, DG13, JL13, CW17, Har21}. In particular, Farina-Hasegawa \cite{FH19} proved Liouville type results for stable solutions to \eqref{1.10} with $s=1$ and a larger class of weights $h$ which cover many existing results.

\medskip Davila-Dupaigne-Wei \cite{DDW17} examined the equation \eqref{1.10} with $h\equiv 1$, and classified finite Morse index solutions in the autonomous case. The approach developed in \cite{DDW17} is based on the monotonicity formula and some energy estimates, this approach was adapted to classify stable or stable solutions outside a compact set to \eqref{1.10} with $h(x)= |x|^\ell$, $\ell\geq 0$ (see \cite{FW16,FW17}) as well as for fractional elliptic equation involving advection term (see \cite{RC19}).

\medskip
Here we obtain classification result for the frational equation by the study of the system. In fact, when $p=q,$ using Souplet type estimate \eqref{2.7}, the system \eqref{1.1} is reduced to the following fractional Lane Emden equation
\begin{align}\label{1.12}
 (-\Delta)^s u = h(x)u^p, \quad u>0 \quad \mbox{in }\;\mathbb{R}^N.
\end{align}
 As a consequence of Theorem \ref{main1}, we can claim
\begin{cor}\label{main2} Suppose that $h$ satisfies $(H)$ and let $p>1.$
\begin{enumerate}
 \item
 If $\frac{4}{3}< p$ then \eqref{1.12} has no stable solution if
\begin{align}\label{1.13}
N<2s+\frac{2(2s+\ell)}{p-1}\left(p+\sqrt{p^2-p}\right).
 \end{align}
 In particular, if $N \leq 10s+4\ell,$ then \eqref{1.12} has no stable
solution for all $\frac{4}{3}< p$.
\item If $1<p\leq \frac{4}{3} $,
\eqref{1.12} has no bounded stable solution for $N$ verifying \eqref{1.13}.
\end{enumerate}
Therefore, there is no bounded stable solution of \eqref{1.12} for all $p > 1$ if $N \leq 10s+4\ell.$
\end{cor}

\medskip
In \cite{TuanHoang21}, Duong-Nguyen studied a more general fractional equation, obtaining similar results of \cite[Theorem 1.5]{DF10} for the Laplacian case.
\begin{equation}\label{1.14}
(-\Delta)^su=f(u)\mbox{ in } \R^N, \quad 0<s<1,
\end{equation}
where $f\in C^0(\mathbb R_+)\cap C^2(\mathbb R_+^*)$ and for  $t>0$, let
$$q(t)= \dfrac{f'^2}{ff''}(t), \;\; \mbox{ if }ff''(t)\not=0; \quad q(t) =
+\infty,\;\;\mbox{ if }ff''(t)=0.$$
Duong-Nguyen proved that if moreover $f$ is nondecreasing, convex, $f>0$ in $\mathbb R_+^*$ and
$q_0= \lim_{t\to 0^+} q(t)\in [1,+\infty]$ exists. Then \eqref{1.14} has no nontrivial bounded nonnegative stable solution if one of the following conditions is satisfied: $(i)$ $N<10s$; $(ii)$ $N=10s$ and $q_0 > 1$; $(iii)$ $N>10s$ and $p_0$ the conjugate exponent of $q_0$ satisfies $p_0<p_c(N,s),$ where
$$p_c(N,s)=\frac{(N-2s)^2-4sN+8s\sqrt{s(N-s)}}{(N-10s)(N-2s)}.$$

\begin{itemize}
\item If $f(u)=u^p$, $p\geq 1$ in \eqref{1.14}, then $p=p_0.$ We can check that the range of $p$ for the nonexistence of stable solution provided by Corollary \ref{main2} with $h\equiv 1$, is the same given by the above result of Duong-Nguyen.
\item In Corollary \ref{main2}, we prove the classification result for stable solution of \eqref{1.12} for $p>\frac43,$ without assuming the boundedness of $u$.
\item Our approach permits to establish Liouville type results for stable solutions to \eqref{1.12} with weights that are not covered (up to our knowledge) by previous works.
\end{itemize}

 This paper is organized as follows. In Section \ref{s3}, we prove comparison properties between $u$ and $v$ of solutions to \eqref{1.1}, and integral estimates derived from the stability. The proof of Theorem \ref{main1} and Corollary \ref{main2} are given in Section \ref{s4}. In the following, $C$ will denote a generic positive constant independent on $(u,v)$, which could be changed from one line to another. The ball of center $0$ and radius $r > 0$ will be denoted by $B_r$.

\section{Preliminaries}\label{s3}
\setcounter{equation}{0}
In this section, we introduce some preliminary results for solutions to the system \eqref{1.1},  as integral estimates; comparison property of $u$, $v$; and an integral inequality derived from the stability.

\medskip
We shall use a standard tool due to Caffarelli-Silvestre \cite{CS07} which transforms the nonlocal system \eqref{1.1} to a degenerate but local elliptic system, with nonlinear Neumann boundary condition in the half space $\mathbb{R}^{N+1}_+.$ More precisely, let $(u,v)$ be a solution of \eqref{1.1}, the extension $(U,V)$ of $(u,v)$  in the sense of \cite{CS07} is defined as follows:
for  $(x,t)\in\mathbb R^{N+1}_+,$
\begin{equation}\label{2.1}
U(x,t)=\int_{\mathbb R^N}P_s(x-z,t)u(z)dz,\quad V(x,t)=\int_{\mathbb R^N}P_s(x-z,t)v(z)dz
\end{equation}
where $P_s(x,t)$ is the Poisson kernel
$$P_s(x,t)=C(N,s)\frac{t^{2s}}{(|x|^2+t^2)^{\frac{N+2s}{2}}}$$
and $C(N,s)$ is a normalization constant.
Then $U, V \in C^2(\mathbb R^{N+1}_+)\cap C(\overline{\mathbb R^{N+1}_+})$, $t^{1-2s}\partial_t U, t^{1-2s}\partial_t V \in C(\overline{\mathbb R^{N+1}_+})$ satisfy
\begin{align}\label{2.2}
\begin{cases}
-{\rm div}(t^{1-2s}\nabla U)=0&\mbox{ in }\mathbb R^{N+1}_+\\
U=u&\mbox{ on }\partial\mathbb R^{N+1}_+\\
-\lim\limits_{t\to 0}t^{1-2s}\partial _t U=\kappa_s(-\Delta)^su&\mbox{ on }\partial\mathbb R^{N+1}_+
\end{cases}
\end{align}
and
\begin{align}\label{2.3}
\begin{cases}
-{\rm div}(t^{1-2s}\nabla V)=0&\mbox{ in }\mathbb R^{N+1}_+\\
V=v&\mbox{ on }\partial\mathbb R^{N+1}_+\\
-\lim\limits_{t\to 0}t^{1-2s}\partial _t V=\kappa_s(-\Delta)^sv&\mbox{ on }\partial\mathbb R^{N+1}_+.
\end{cases}
\end{align}
Here $\kappa_s=\frac{\Gamma(1-s)}{2^{2s-1}\Gamma(s)}$  and $\Gamma$ is the usual Gamma function. For any $W\in C(\overline{R^{N+1}})$ and $r > 0$, we define
$$\overline{W}(r):=\frac{1}{r^{N+1-2s}}\displaystyle{\int_{\partial^+ B_r}}y^{1-2s}W,$$
where $B_r^+:=B_r \cap \{y>0\}$ is the half-ball with spherical part of boundary $\partial^+B_r= \partial B_r\cap \{y>0\}.$

\medskip
Inspired by \cite{YZ19_b, SZ96, MPO}, we establish the following \textit{a priori} integral estimates for solutions to the fractional Lane-Emden system \eqref{1.1}.

\begin{lemma}\label{l.2.1}
Let $p, q \geq1$, $pq > 1$. Suppose that $h$ satisfies $(H).$ Then, there exists a positive constant $C$ depending only on $N,s,p$ and $q$ such that for any solution $(u,v)$ of \eqref{1.1} and any $R > 0$, there hold
\begin{align}\label{2.4}
\begin{split}
\int_{B_R} h(x) u^q(x) dx \leq CR^{N- \frac{2s q(p+1)}{pq-1}-\frac{ \ell(q+1)}{pq-1}},\quad \int_{B_R} h(x) v^p(x) dx \leq CR^{N- \frac{2s p(q+1)}{pq-1}-\frac{ \ell(p+1)}{pq-1}}.
\end{split}
\end{align}
\end{lemma}

\medskip
\noindent {\textbf{Proof.}
 Let $\varphi_0 \in C_c^\infty(B_2)$ be a cut-off function verifying $0 \leq \varphi_0 \leq 1$, and $\varphi_0=1$ for $x  \in B_1$. Let $R > 0$ and consider
  $\psi := \varphi_0(R^{-1}x).$ By \cite[Lemma 2.1]{YZ19_b}, we have for $m\geq 1,$
  \begin{align}\label{2.5}
   (-\D)^s \psi ^{m}\leq m R^{-2s}[\psi ^{m-1}(-\D)^s \psi](R^{-1}x).
  \end{align}
    Multiplying the equation $(-\D)^s u=h(x)v^{p}$ by $\psi ^{m}$ and integrating by parts,
there holds
\begin{align*}
   \int_{\R^N}h(x)v^{p}\psi^{m} dx = \int_{\R^N}u(-\Delta)^s\psi^{m} dx \leq \frac{C}{R^{2s}}\int_{B_{2R}} u \psi^{m-1} dx.
\end{align*}
By H\"older's inequality, we obtain
\begin{align*}
  \int_{\R^N}h(x)v^{p}\psi^{m}dx  \leq& \; \frac{C}{R^{2s}}\left(\int_{B_{2R}} h(x)^{-\frac{q'}{q}}dx\right)^{\frac{1}{q'}}\left(\int_{B_{2R}}h(x) u^{q} \psi^{(m-1)q}dx\right)^{\frac{1}{q}},
\end{align*}
where $\frac{1}{q}+\frac{1}{q'}=1.$ From $(H)$ we deduce that
\begin{align*}
  \int_{\R^N}h(x) v^{p}\psi^{m}dx  \leq& \; C R^{\frac{N}{q'} -\frac{\ell}{q}-2s}\left(\int_{B_{2R}}h(x) u^{q} \psi^{(m-1)q}dx\right)^{\frac{1}{q}}.
\end{align*}
 Similarly, using $(-\Delta)^s v = h(x) u^q$, for $k \geq 1$,
\begin{align*}
  \int_{\R^N}h(x) u^{q}\psi^{k} dx \leq& \; C R^{\frac{N}{p'}-\frac{\ell}{p} -2s}\left(\int_{B_{2R}}h(x) v^{p} \psi^{(k-1)p}dx\right)^{\frac{1}{p}},
\end{align*}
where $\frac{1}{p}+\frac{1}{p'}=1.$ Take now $k$ and $m$ large verifying $ m \leq (k-1)p$ and $k \leq (m-1)q.$ Combining the two
above inequalities, we get
\begin{align*}
  \int_{\R^N}h(x)v^{p}\psi^{m}dx & \leq C R^{\frac{N}{q'} -\frac{\ell}{q}-2s} R^{\big(\frac{N}{p'} -\frac{\ell}{p}-2s\big)\frac{1}{q}}\left(\int_{B_{2R}}h(x) v^{p} \psi^{(k-1)p}dx\right)^{\frac{1}{pq}}\\
  & \leq C R^{N - \frac{N}{pq} -\frac{\ell(p+1)}{pq}- \frac{2s(q+1)}{q}}\left(\int_{\R^N}h(x) v^{p} \psi^m dx\right)^{\frac{1}{pq}}.
\end{align*}
Hence
\begin{align*}
  \int_{B_R}h(x)v^{p} dx\leq \int_{\R^N}h(x)v^{p}\psi^{m}dx \leq C R^{N-\frac{2s(q+1)p}{pq-1}-\frac{(p+1)\ell}{pq-1}}.
\end{align*}
Similarly, we obtain the estimate for $u$. \qed

\medskip
An immediate consequence of Lemma \ref{l.2.1} is the following non-existence theorem for the fractional weighted Lane Emden system \eqref{1.1}.
\begin{cor}\label{main2.1}
Let $p, q \geq1$, $pq > 1$. Suppose that $h$ satisfies $(H).$ Then, there exists no solution to \eqref{1.1} if
$N< \frac{2spq + \ell}{pq-1}+ \frac{2s +\ell}{pq-1}\times\max(p, q).$
\end{cor}

Using Lemma \ref{l.2.1}, we can also give an alternative proof of \cite[Proposition 3.2]{TuanHoang21} under the same assumptions of Lemma \ref{l.2.1}. Consider for example
\begin{align}\label{2.6}
\begin{split}
& \int_{\R^N} h(x) u^q(x)\rho(\frac{x}{R}) dx \leq CR^{N- \frac{2s q(p+1)}{pq-1}-\frac{ \ell(q+1)}{pq-1}}.
\end{split}
\end{align}
Here and after $\rho(x) := (1 + |x|^2)^{-\frac{N+2s}{2}}.$ Indeed, let $k\in \mathbb{N}^{*}.$ Then $\rho(\frac{x}{R})\leq 1$ for any $x\in B_R$ and $\rho(\frac{x}{R})\leq 2^{-(k-1)(N+2s)}$ for any $x\in B_{2^{k}R}\setminus B_{2^{k-1}R}.$ Applying Lemma \ref{l.2.1}, it follows that
\begin{align*}
   \int_{\R^N} h(x) u^q(x)\rho(\frac{x}{R})dx  &= \int_{B_R} h(x) u^q(x)\rho(\frac{x}{R})dx+ \displaystyle{\sum_{k\in \mathbb{N}^{*}}} \int_{B_{2^{k}R}\setminus B_{2^{k-1}R}} h(x) u^q(x)\rho(\frac{x}{R})dx  \\
   &\leq \int_{B_R} h(x) u^q(x)dx + \displaystyle{\sum_{k\in \mathbb{N}^{*}}} 2^{-(k-1)(N+2s)}\int_{B_{2^{k}R}\setminus B_{2^{k-1}R}} h(x) u^q(x)dx\\
   &\leq CR^{N- \frac{2s q(p+1)}{pq-1}-\frac{ \ell(q+1)}{pq-1}}\left[1+\displaystyle{\sum_{k\in \mathbb{N}^{*}}} 2^{-(k-1)(N+2s)+k(N- \frac{2s q(p+1)}{pq-1}-\frac{ \ell(q+1)}{pq-1})}\right]\\
   &\leq CR^{N- \frac{2s q(p+1)}{pq-1}-\frac{ \ell(q+1)}{pq-1}}.
\end{align*}

The following is a comparison result between the components $u$, $v$ of solutions to the system \eqref{1.1}.
\begin{proposition}\label{p.2.1}
 Let $p\geq q \geq 1$ and $pq > 1$. Suppose that $h$ satisfies $(H)$ and $(u,v)$ is a solution to \eqref{1.1}. There holds then
\begin{equation}\label{2.7}
	v^{p+1}\leq \frac{p+1}{q+1}u^{q+1}.
	\end{equation}
If moreover $v$ is bounded, then
  \begin{align}\label{2.8}
u\leq \|v\|_\infty^\frac{p-q}{q+1}v.
 \end{align}
\end{proposition}

\medskip\noindent
{\textbf{Proof.}
The proof adapt an idea of \cite{YZ19_b}, originally coming from \cite{QS12}. Let $w:=v-\l u^{\sigma},$ where $\sigma=\frac{q+1}{p+1}$ and $\l =\sigma ^{-\frac{1}{p+1}},$ the proof of \eqref{2.7} consists to show that
\begin{align}\label{2.9}
 (-\Delta)^sw\leq 0 \quad \mbox{in the set}\quad \{w\geq 0\}.
\end{align}
Indeed, let $W$ be the extension of $w$ in the sense of \eqref{2.2}-\eqref{2.3}. Using \eqref{2.9}, then $W$ satisfies
\begin{align}\label{2.10}
\begin{cases}
-{\rm div}(t^{1-2s}\nabla W)=0&\mbox{ in }\mathbb R^{N+1}_+\\
-\lim\limits_{t\to 0}t^{1-2s}\partial _t W=\kappa_s(-\Delta)^sw\leq 0&\mbox{ on }\{W\geq 0\}\bigcap\partial\mathbb R^{N+1}_+.
\end{cases}
\end{align}
Moreover, by the integral estimate \eqref{2.4}, there exist $r_i\rightarrow +\infty,$ $x_i\in B_{r_i}$ such that
$$v^p(x_i)r_i^N \leq Cr_i^{{N- \frac{2s p(q+1)}{pq-1}-\frac{ \ell(p+1)}{pq-1}}},$$ which implies $\lim_{i\rightarrow +\infty}v(x_i)=0.$ Applying \cite[Lemma 3.2]{YZ19_b}, there hold $0\leq \lim_{r \rightarrow +\infty}\overline{V}(r)\leq Av(x_i)$, $\forall\; i \in \N$ for some positive constant $A.$ Tending $i$ to $+\infty,$ we obtain $\lim_{r \rightarrow +\infty}\overline{V}(r)=0.$ Furthermore, $0\leq \overline{W_+}(r)\leq \overline{V}(r),$ where $W_+=\max(W,0).$ Hence, $\lim_{r \rightarrow +\infty}\overline{W_+}(r)=0.$ Following the same lines in the proof of \cite[Lemma 3.1]{YZ19_b}, we derive that $W\leq 0$ and hence $w\leq0,$ i.e. \eqref{2.7}.

\medskip
To get \eqref{2.9}, consider the concave function $t^\sigma$ in $\R_+$,
$$u^{\sigma}(x)-u^{\sigma}(y)\geq \sigma u^{\sigma-1}(x)(u(x)-u(y)), \quad \forall\; x, y.$$
Hence
	$$(-\Delta)^su^{\sigma}(x)=\int_{\mathbb R^N}\frac{u^{\sigma}(x)-u^{\sigma}(y)}{|x-y|^{N+2s}}dy\geq \sigma u^{\sigma-1}(x)(-\Delta)^su(x).$$
It follows that
\begin{align*}
 (-\Delta)^sw  =(-\Delta)^sv-\l(-\Delta)^s u^{\sigma} \leq h u^q-\l h\sigma u^{\sigma-1}v^p  & = h\left[u^q- \l^{-p} u^{\sigma-1}v^p\right] \notag\\
 & = \l^{-p} u^{\sigma-1}h \left[(\l u^{\sigma})^p-v^p\right]\notag\\
   &\leq 0\;\mbox{on the set }\{w\geq 0\}.
\end{align*}
So we are done. To prove \eqref{2.8}, consider $w= u- \l v$ with $\l = \|v\|_{\infty}^\frac{p-q}{q+1}$ and we will establish again \eqref{2.9}. As $p \geq q$ and $v$ is bounded, there holds
\begin{align*}
(-\Delta)^s w = h(x) v^p- \l h(x)u^q \leq h(x)\left(v^p-\l u^q\right) & =h(x)\left[ \left(\frac{v}{ \|v\|_{\infty}}\right)^{p} \|v\|_{\infty}^{p} -\l u^q \right]\\
 &\leq h(x)\left[ \left(\frac{v}{ \|v\|_{\infty}}\right)^{q} \|v\|_{\infty}^{p} -\l u^q \right]\\
&=h(x) \|v\|_{\infty}^{p - q} \left(v^q-\frac{ \l u^q}{\|v\|_{\infty}^{p - q}}\right)\\
&= h(x)\|v\|_{\infty}^{p - q} \left(v^q-\l^{-q}u^q \right).
\end{align*}
Therefore, we get \eqref{2.9} and the proof is completed.  \qed
\begin{rem}\label{r.2.1} Let $(u,v)$ be a solution of \eqref{1.1} and $(U,V)$ be the extension of $(u,v)$ in the sense of \eqref{2.2}--\eqref{2.3}. With the assumptions of Proposition \ref{p.2.1}, there holds
\begin{equation}\label{2.11}
	\frac{V^{p+1}}{p+1}\leq \frac{U^{q+1}}{q+1}.
	\end{equation}
Indeed, using the same notations as above, we have $\frac{1}{\sigma}\geq 1.$ Hence, by Jensen's inequality and \eqref{2.7}, we get for $(x,t)\in\mathbb R^{N+1}_+,$
\begin{align*}
  (V(x,t))^{\frac{1}{\sigma}}&\leq \int_{\mathbb R^N}P_s(x-z,t)(v(z))^{\frac{1}{\sigma}}dz \leq l^{\frac{1}{\sigma}}\int_{\mathbb R^N}P_s(x-z,t)u(z)dz= l^{\frac{1}{\sigma}}U(x,t).
\end{align*}
\end{rem}

At last, using the Definition \ref{d.1.1} of stability, we can derive the following estimate which is crucial for our analysis. Its proof comes from ideas in \cite{CG14, DFS13, Cow13, TuanHoang21} and is very similar to the mentioned works, so we omit the details.
\begin{lemma}\label{l.2.2}
	Let  $(u,v)$ be a stable solution of \eqref{1.1}. Then for all $\phi\in C_c^\infty(\mathbb R^N)$, we have
\begin{align}\label{2.12}
  \sqrt{pq}\int_{\mathbb R^N}h(x) u(x)^{\frac{q-1}{2}}v^{\frac{p-1}{2}}\phi(x)^2dx\leq \frac{c_{N,s}}{2}\int_{\mathbb R^N}\int_{\mathbb R^N}\frac{(\phi(x)-\phi(y))^2}{|x-y|^{N+2s}}dydx.
\end{align}
\end{lemma}

\section{Proof of Theorem \ref{main1} and Corollary \ref{main2}.}\label{s4}
 Assume that $(u,v)$ is a stable solution of \eqref{1.1} and $h$ satisfies $(H).$
Denote by $U$, $V$ the extension of $u$ and $v$ in the sense of \eqref{2.2}--\eqref{2.3} and define $\zeta(x,\,t):= (1+|x|^2+t^2)^{-\frac{N+2s}{4}}$ an extension of $\rho(x)^{\frac12} $ on $\mathbb{R}_+^{N+1}.$
\begin{lemma}\label{l.3.1}
For any $\gamma>\frac{q+1}{2}$ satisfying $L(\gamma)<0$ and $\Phi\in C_c^\infty(\mathbb R^{N+1}_+)$, there exists $C>0$ such that
\begin{align}\label{3.1}
\int_{\mathbb R^{N+1}_+}|\nabla (U^{\frac{\gamma}{2}} \zeta\Phi)|^2t^{1-2s}dxdt\leq C\int_{\mathbb R^{N+1}_+}U^{\gamma}|\nabla (\zeta\Phi) |^2 t^{1-2s}dxdt
\end{align}
where \begin{align}\label{3.2}
 L(\gamma):=\gamma^4-16\frac{pq(q+1)}{p+1}\gamma^2+16\frac{pq(q+1)(p+q+2)}{(p+1)^2}\gamma-
 16\frac{pq(q+1)^2}{(p+1)^2}.
\end{align}
\end{lemma}

\medskip
\noindent
{\textbf{Proof.}
	Let $\Phi\in C^\infty_c(\mathbb R^{N+1}_+)$ be a test function and define $\phi(x)=\Phi(x,0)\in C^\infty_c(\mathbb R^N)$. Let $\gamma >1.$ Multiplying the first equation in \eqref{2.2} by $U^{\gamma-1}(\zeta\Phi)^2$ and integrating by parts, we get
	\begin{align}\label{3.3}
	\begin{split}
	&\quad \kappa_s\int_{\mathbb R^N}h(x)v(x)^pu(x)^{\gamma-1}\rho(x)\phi(x)^2dx\\
& =\int_{\mathbb R^{N+1}_+}\nabla U\cdot\nabla(U^{\gamma-1}(\zeta\Phi)^2)t^{1-2s}dxdt\\
&=(\gamma-1)\int_{\mathbb R^{N+1}_+}|\nabla U|^2U^{\gamma-2}(\zeta\Phi)^2t^{1-2s}dxdt+\frac{4}{\gamma}\int_{\mathbb R^{N+1}_+}\zeta\Phi\nabla(U^{\frac{\gamma}{2}})\cdot\nabla (\zeta\Phi) U^{\frac{\gamma}{2}} t^{1-2s}dxdt\\
&=\frac{4(\gamma-1)}{\gamma^2}\int_{\mathbb R^{N+1}_+}|\nabla U^{\frac{\gamma}{2}}|^2(\zeta\Phi)^2t^{1-2s}dxdt+\frac{4}{\gamma}\int_{\mathbb R^{N+1}_+}\zeta\Phi\nabla(U^{\frac{\gamma}{2}})\cdot\nabla (\zeta\Phi) U^{\frac{\gamma}{2}} t^{1-2s}dxdt.
	\end{split}
	\end{align}
Furthermore,  there holds
\begin{align}\label{3.4}
\begin{split}
\int_{\mathbb R^{N+1}_+}|\nabla U^{\frac{\gamma}{2}}|^2(\zeta\Phi)^2t^{1-2s}dxdt=&\int_{\mathbb R^{N+1}_+}|\nabla (U^{\frac{\gamma}{2}}\zeta\Phi)|^2t^{1-2s}dxdt\\
&-2\int_{\mathbb R^{N+1}_+}\zeta\Phi\nabla(U^{\frac{\gamma}{2}})\cdot\nabla (\zeta\Phi) U^{\frac{\gamma}{2}} t^{1-2s}dxdt\\
&-\int_{\mathbb R^{N+1}_+}U^{\gamma}|\nabla (\zeta\Phi) |^2 t^{1-2s}dxdt.
\end{split}
\end{align}
Hence, using the Cauchy-Schwarz inequality, for any $\epsilon >0,$ we have
	\begin{align}\label{3.5}
\begin{split}
\kappa_s\int_{\mathbb R^N}h(x)v(x)^pu(x)^{\gamma-1}\rho(x)\phi(x)^2dx
=&\frac{4(\gamma-1)}{\gamma^2}\int_{\mathbb R^{N+1}_+}|\nabla (U^{\frac{\gamma}{2}}\zeta\Phi)|^2t^{1-2s}dxdt\\
&-\frac{4(\gamma-2)}{\gamma^2}\int_{\mathbb R^{N+1}_+}\zeta\Phi\nabla(U^{\frac{\gamma}{2}})\cdot\nabla (\zeta\Phi) U^{\frac{\gamma}{2}} t^{1-2s}dxdt\\
&-\frac{4(\gamma-1)}{\gamma^2}\int_{\mathbb R^{N+1}_+}U^{\gamma}|\nabla (\zeta\Phi) |^2 t^{1-2s}dxdt\\
=&\frac{4(\gamma-1)}{\gamma^2}\int_{\mathbb R^{N+1}_+}|\nabla (U^{\frac{\gamma}{2}}\zeta\Phi)|^2t^{1-2s}dxdt\\
&-\frac{4(\gamma-2)}{\gamma^2}\int_{\mathbb R^{N+1}_+}\nabla(\zeta\Phi U^{\frac{\gamma}{2}})\cdot\nabla (\zeta\Phi) U^{\frac{\gamma}{2}} t^{1-2s}dxdt\\
&-\frac{4}{\gamma^2}\int_{\mathbb R^{N+1}_+}U^{\gamma}|\nabla (\zeta\Phi) |^2 t^{1-2s}dxdt\\
\geq& \frac{4(\gamma-1)}{\gamma^2}(1-\epsilon)
\int_{\mathbb R^{N+1}_+}|\nabla (U^{\frac{\gamma}{2}}\zeta\Phi)|^2t^{1-2s}dxdt\\
&-\left(\frac{4}{\gamma^2}+\frac{C_\gamma}{\epsilon}\right)\int_{\mathbb R^{N+1}_+}U^{\gamma}|\nabla (\zeta\Phi) |^2 t^{1-2s}dxdt,
\end{split}
\end{align}
Denote by $A(\gamma, \epsilon):= \frac{4\sqrt{pq}(\gamma-1)}{\gamma^2}(1-\epsilon).$ It follows that
\begin{align}\label{3.6}
\begin{split}
\frac{1}{\sqrt{pq}}\int_{\mathbb R^{N+1}_+}|\nabla (U^{\frac{\gamma}{2}}\zeta\Phi)|^2t^{1-2s}dxdt
\leq& \frac{\kappa_s }{A(\gamma, \epsilon)}\int_{\mathbb R^N}h(x)v(x)^pu(x)^{\gamma-1}\rho(x)\phi(x)^2dx\\
&+C_{\gamma, \epsilon}\int_{\mathbb R^{N+1}_+}U^{\gamma}|\nabla (\zeta\Phi) |^2 t^{1-2s}dxdt.
\end{split}
\end{align}
Similarly, multiplying the first equation in \eqref{2.3} by $V^{\gamma-1}(\zeta\Phi)^2$ and integrating by parts, we obtain
\begin{align}\label{3.7}
\begin{split}
\frac{1}{\sqrt{pq}}\int_{\mathbb R^{N+1}_+}|\nabla (V^{\frac{\gamma}{2}}\zeta\Phi)|^2t^{1-2s}dxdt
\leq& \frac{\kappa_s}{A(\gamma, \epsilon)}\int_{\mathbb R^N}h(x) u(x)^qv(x)^{\gamma-1}\rho(x)\phi(x)^2dx \\
&+C_{\gamma, \epsilon}\int_{\mathbb R^{N+1}_+}V^{\gamma}|\nabla (\zeta \Phi) |^2 t^{1-2s}dxdt.
\end{split}
\end{align}
Combining \eqref{3.6} and \eqref{3.7}, we derive that, for any $\gamma_1, \gamma_2>1,$
 \begin{align}
 \label{3.8}
  \begin{split} &A(\gamma_1, \epsilon)^{\frac{2\gamma_1}{q+1}} I_1+ I_2\\
 := & \; A(\gamma_1, \epsilon)^{\frac{2\gamma_1}{q+1}} \frac{1}{\sqrt{pq}}\int_{\mathbb R^{N+1}_+}|\nabla (U^{\frac{\gamma_1}{2}}\zeta\Phi)|^2t^{1-2s}dxdt+
 \frac{1}{\sqrt{pq}} \int_{\mathbb R^{N+1}_+}|\nabla (V^{\frac{\gamma_2}{2}}\zeta\Phi)|^2t^{1-2s}dxdt\\
  \leq & \;\kappa_s A(\gamma_1, \epsilon)^\frac{2\gamma_1-1-q}{q+1}\int_{\mathbb R^N}h(x)v(x)^pu(x)^{\gamma_1-1}\rho(x)\phi(x)^2dx\\
&+\frac{\kappa_s}{A(\gamma_2, \epsilon)} \int_{\mathbb R^N}h(x)u(x)^qv(x)^{\gamma_2-1}\rho(x)\phi(x)^2dx\\
&+C_{\epsilon}\int_{\mathbb R^{N+1}_+}\left(U^{\gamma_1}+V^{\gamma_2}\right)|\nabla (\zeta\Phi) |^2 t^{1-2s}dxdt.
\end{split}
 \end{align}
Fix now
\begin{equation}\label{3.9}
\gamma_2=\frac{(p+1)}{q+1}\gamma_1 \quad \iff \quad
 \gamma_2-1=\frac{p+1}{q+1}(\gamma_1-1)+\frac{p-q}{q+1}.
\end{equation}
Let $\gamma_1> \frac{q+1}{2},$ by Young's inequality, there holds
\begin{align*}
& \frac{\kappa_s}{A(\gamma_2, \epsilon)} \int_{\mathbb R^N}h(x) u(x)^qv(x)^{\gamma_2-1}\rho(x)\phi(x)^2dx\\
= & \; \frac{\kappa_s}{A(\gamma_2,\epsilon)}\int_{\mathbb{R}^N}h(x) u(x)^{\frac{q-1}{2}}v(x)^{\frac{p-1}{2}}v(x)^{\frac{(p+1)\gamma_1}{q+1}-
\frac{p+1}{q+1}\left(\frac{q+1}{2}\right)}u(x)^{\frac{q+1}{2}}\rho(x)\phi(x)^2 dx \\
= & \;\frac{\kappa_s}{A(\gamma_2,\epsilon)}\int_{\mathbb{R}^N}h(x) u(x)^{\frac{q-1}{2}}v(x)^{\frac{p-1}{2}}v(x)^{\gamma_2\frac{2\gamma_1-q-1}{2\gamma_1}}u(x)^{\frac{q+1}{2}}\rho(x)\phi(x)^2dx\\
\leq & \; \frac{2\gamma_1-q-1}{2\gamma_1}\kappa_s\int_{\mathbb{R}^N}h(x) u(x)^{\frac{q-1}{2}}v(x)^{\frac{p-1}{2}}v(x)^{\gamma_2}\rho(x)\phi(x)^2dx\\
 &+\frac{q+1}{2\gamma_1}A(\gamma_2,\epsilon)^{-\frac{2\gamma_1}{q+1}}\kappa_s\int_{\mathbb{R}^N}h(x) u(x)^{\frac{q-1}{2}}v(x)^{\frac{p-1}{2}}u(x)^{\gamma_1}\rho(x)\phi(x)^2dx.
\end{align*}
Choosing the test function $u^\frac{\gamma_1}{2}\rho_{N+2s}^{\frac12}\phi$ (resp. $v^\frac{\gamma_2}{2}\rho_{N+2s}^{\frac12}\phi$)  in the stability inequality \eqref{2.12} and using the fact that $U^{\frac{\gamma_1}{2}}\zeta\Phi$ (resp. $V^{\frac{\gamma_2}{2}}\zeta\Phi$) has the trace $u^{\frac{\gamma_1}{2}}\rho_{N+2s}^{\frac12}\phi$ (resp. $v^{\frac{\gamma_2}{2}}\rho_{N+2s}^{\frac12}\phi$) on $ \partial \R^{N+1}_+$, one gets
\begin{align}\label{3.10}
\begin{split}
 \kappa_s\sqrt{pq}\int_{\mathbb R^N}h(x) v(x)^\frac{p-1}{2}u(x)^{\frac{q-1}{2}}u(x)^{\gamma_1}\rho(x)\phi(x)^2dx
&\leq \kappa_s\|u^\frac{\gamma_1}{2}\rho^{\frac12}\phi\|_{\overset{.}{H}^s(\mathbb{R}^N)}\\
&\leq \int_{\mathbb R^{N+1}_+}|\nabla (U^{\frac{\gamma_1}{2}}\zeta\Phi)|^2t^{1-2s}dx=\sqrt{pq}I_1,
\end{split}
\end{align}
and
\begin{align}\label{3.11}
\begin{split}
 \kappa_s\sqrt{pq}\int_{\mathbb R^N}h(x) v(x)^\frac{p-1}{2}u(x)^{\frac{q-1}{2}}v(x)^{\gamma_2}\rho(x)\phi(x)^2dx
&\leq \kappa_s\|v^\frac{\gamma_2}{2}\rho_{N+2s}^{\frac12}\phi\|_{\overset{.}{H}^s(\mathbb{R}^N)}\\
&\leq \int_{\mathbb R^{N+1}_+}|\nabla (V^{\frac{\gamma_2}{2}}\zeta\Phi)|^2t^{1-2s}dxdt=\sqrt{pq}I_2.
\end{split}
\end{align}
Hence,
\begin{align*}
 \frac{\kappa_s}{A(\gamma_2,\epsilon)} \int_{\mathbb R^N}h(x) u(x)^qv(x)^{\gamma_2-1}\rho(x)\phi(x)^2dx
\leq  \frac{2\gamma_1-q-1}{2\gamma_1}I_2
 +\frac{q+1}{2\gamma_1}\left(A(\gamma_2,\epsilon)\right)^{-\frac{2\gamma_1}{q+1}}I_1.
\end{align*}
Similarly, we can prove that
 \begin{equation}\nonumber
  \kappa_s A(\gamma_1, \epsilon)^\frac{2\gamma_1-1-q}{q+1}\int_{\mathbb{R}^N}h(x) v(x)^p
u(x)^{\gamma_1-1}\rho(x)\phi^2dx \leq \frac{q+1}{2\gamma_1} I_2
  + \frac{2\gamma_1-q-1}{2\gamma_1}A(\gamma_1, \epsilon)^{\frac{2\gamma_1}{q+1}}I_1.
   \end{equation}
   Combining the above two estimates with \eqref{3.8}, we derive that
\begin{align*}
\begin{split}
A(\gamma_1, \epsilon)^{\frac{2\gamma_1}{q+1}}I_1\leq&
\left[\frac{q+1}{2\gamma_1}A(\gamma_2,\epsilon)^{-\frac{2\gamma_1}{q+1}}+\frac{2\gamma_1-q-1}{2\gamma_1}A(\gamma_1, \epsilon)^{\frac{2\gamma_1}{q+1}}\right]I_1 \\
&+ C_{\epsilon}\int_{\mathbb R^{N+1}_+}\left(U^{\gamma_1}+V^{\gamma_2}\right)|\nabla (\zeta\Phi) |^2 t^{1-2s}dxdt,
\end{split}
\end{align*}
hence
\begin{equation*}
\frac{q+1}{2\gamma_1}\left[\left(A(\gamma_1,\epsilon)A(\gamma_2,\epsilon)\right)^{\frac{2\gamma_1}{q+1}}-1\right]
I_1\leq A(\gamma_2,\epsilon)^{\frac{2\gamma_1}{q+1}}C_{\epsilon}\int_{\mathbb R^{N+1}_+}\left(U^{\gamma_1}+V^{\gamma_2}\right)|\nabla (\zeta\Phi) |^2 t^{1-2s}dxdt.
\end{equation*}
Denote  $A_i= A(\gamma_i,0),\,i=1,\,2.$ Suppose that $A_1A_2>1,$ we can choose $\epsilon$ small enough such that $A(\gamma_1,\epsilon)A(\gamma_2,\epsilon)>1,$ and so we obtain
\begin{equation*}
    I_1\leq C\int_{\mathbb R^{N+1}_+}\left(U^{\gamma_1}+V^{\gamma_2}\right)|\nabla (\zeta\Phi) |^2 t^{1-2s}dxdt.
\end{equation*}
On the other hand, by Remark \ref{r.2.1}, there holds $V^{\gamma_2}\leq C U^{\gamma_1}.$ Denoting $\gamma:=\gamma_1,$ we conclude that if $A_1A_2> 1$ and $\gamma > \frac{q+1}{2}$,
\begin{equation*}
 \int_{\mathbb R^{N+1}_+}|\nabla (U^{\frac{\gamma}{2}}\zeta\Phi)|^2t^{1-2s}dxdt  \leq
 C\int_{\mathbb R^{N+1}_+}U^{\gamma}|\nabla (\zeta\Phi) |^2 t^{1-2s}dxdt.
\end{equation*}
Finally, we can check that $A_{1}A_2>1$ is equivalent to $L(\gamma)<0$, the proof is completed. \qed

\medskip
\noindent
\textbf{End of the proof of Theorem \ref{main1}}.
Take $\phi \in C_c^\infty(({-2}, 2))$ satisfying $\phi \equiv 1$ in $[-1,\,1]$. For $(x,\,t)\in \mathbb R^{N+1}_+$ and $R>0,$ we define $\Phi_R\in C_c^\infty(\mathbb R^{N+1}_+)$ by $\Phi_R(x,\,t)=\phi \left(\frac{|(x,\,t)|}{R}\right).$
Let $\gamma_0$ be the largest root of the polynomial $L$ given by \eqref{3.2} and denote by
$$k_s=\frac{N+2-2s}{N-2s}.$$
Fix $\tau>0,$ then there exists a positive integer $m$ satisfying
$$\tau k_s^{m-1}<\gamma_0\leq \tau k_s^{m}.$$
Define $\gamma_1,...,\gamma_m$ as follows
	$$ \gamma_1=\tau k<\gamma_2=\tau kk_s< ...<\gamma_m=\tau kk_s^{m-1}<\gamma_0, $$
	where $k\in [1,k_s]$ will be chosen so that $\gamma_m$ is arbitrarily close to $\gamma_0.$ Suppose that $\tau$ satisfies
$$\tau>\frac{q+1}{2}\;\mbox{ such that }\; L(\gamma)<0 \;\mbox{ for any }\; \gamma \in (\tau,\, \gamma_0),\leqno{(*)}$$
then  $\gamma_1,...,\gamma_m$ satisfy $(*).$ Hence, from \eqref{3.1}  and the Sobolev inequality (see \cite[Proposition 3.1.1]{DMV17}), there holds
	\begin{align*}
\begin{split}
		\left(\int_{B^+_1} U^{\gamma_m k_s}\zeta t^{1-2s}dxdt\right)^{\frac{2}{\gamma_mk_s}}&\leq C \left(\int_{B^+_{2}}U^{\gamma_mk_s}(\zeta\Phi_1)^{2k_s}t^{1-2s}dxdt\right)^{\frac{2}{\gamma_mk_s}}\\
		&\leq 	C\left(\int_{B^+_{2}}|\nabla (U^{\frac{\gamma_m}{2}}(\zeta\Phi_1))|^2t^{1-2s}dxdt\right)^{\frac{2}{\gamma_m}}\\
	&\leq C\left(\int_{B^+_{2}}U^{\gamma_m}|\nabla (\zeta\Phi_1) |^2 t^{1-2s}dxdt\right)^{\frac{2}{\gamma_m}}\\
	&\leq C\left(\int_{B^+_{2}}U^{\gamma_{m-1}k_s} \zeta t^{1-2s}dxdt\right)^{\frac{2}{\gamma_{m-1}k_s}}.
\end{split}
	\end{align*}
 Iterating the above and using H\"older inequality, we arrive at
		\begin{align}\label{3.12}
	\begin{split}
	\left(\int_{B^+_1} U^{\gamma_mk_s} \zeta t^{1-2s}dxdt\right)^{\frac{2}{\gamma_mk_s}}&\leq C\left(\int_{B^+_{2^{m-1}}}U^{\tau k}\zeta t^{1-2s}dxdt\right)^{\frac{2}{\tau k}}\\
	&\leq C\left(\int_{B^+_{2^{m-1}}}U^{\tau k_s}\zeta  t^{1-2s}dxdt\right)^{\frac{2}{\tau k_s}}\\
	&\leq C\left(\int_{B^+_{2^{m}}}U^{\tau }|\nabla (\zeta\Phi_{2^m})|^2 t^{1-2s}dxdt\right)^{\frac{2}{\tau }}\\
&\leq C\left(\int_{\mathbb{R}_+^{N+1}}U^{\tau }|\nabla (\zeta\Phi_{2^m})|^2 t^{1-2s}dxdt\right)^{\frac{2}{\tau }}.
	\end{split}
	\end{align}
Now, we can adapt the proof of \cite[Lemma 2.4]{DuongPham20} (see also \cite[(4.8)--(4.10)]{TuanHoang21} for further explanations) to obtain the following integral estimate from $\mathbb{R}_+^{N+1}$ to $\mathbb{R}^N:$
\begin{align}\label{3.13}
  \int_{\mathbb{R}_+^{N+1}}U^{\tau }|\nabla (\zeta\Phi_{2^m})|^2 t^{1-2s}dxdt\leq \int_{\mathbb{R}^{N}}(u(y))^{\tau}\rho(y) dy\leq \int_{\mathbb{R}^{N}}h(y)(u(y))^{\tau}\rho(y) dy,
\end{align}
where, in the last inequality, we used $h\geq C > 0$ in $\R^N$.
We deduce from \eqref{3.12} that
\begin{align}\label{3.14}
\left(\int_{B^+_1} U^{\gamma_mk_s}\zeta t^{1-2s}dxdt\right)^{\frac{2}{\gamma_mk_s}}&\leq C\left(\int_{\mathbb R^N} h(y)(u(y))^{\tau}\rho(y) dy\right)^{\frac{2}{\tau}}.
\end{align}
Let $R>1$. The functions
$$u_R(x) = R^{\frac{2s(p+1)}{pq-1}+\frac{\ell(p+1)}{pq-1}} u(R x)\mbox{ and } v_R(x) = R^{\frac{2s(q+1)}{pq-1}+\frac{\ell(q+1)}{pq-1}} v(Rx),$$
form a  solution of \eqref{1.1} with $h$ is replaced by $\frac{h(Rx)}{R^\ell}.$ We use a scaling argument, replacing $U(x,t),$ $u(y)$ and $h(y)$ in \eqref{3.14} by $U_R:=R^{\frac{2s(p+1)}{pq-1}+\frac{\ell(q+1)}{pq-1}} U(R x,Rt),$ $u_R(y)$ and $\frac{h(Ry)}{R^\ell},$ we deduce that
\begin{align*}
\left(\int_{B^+_1} U^{\gamma_mk_s}(Rx,Rt)t^{1-2s}dxdt\right)^{\frac{2}{\gamma_mk_s}}&\leq C\left(\int_{\mathbb R^N}\frac{h(Rx)}{R^\ell} (u(Ry))^{\tau}\rho(y) dy\right)^{\frac{2}{\tau}}.
\end{align*}
Hence,
\begin{align}\label{3.15}
\begin{split}
\left(R^{-N-2+2s}\int_{B^+_R} U^{\gamma_mk_s}(x,t)t^{1-2s}dxdt\right)^{\frac{2}{\gamma_mk_s}}&\leq C\left(\int_{\mathbb R^N}\frac{h(Rx)}{R^\ell} (u(Ry))^{\tau}\rho(y) dy\right)^{\frac{2}{\tau}}\\
&\leq C\left(R^{-N-\ell}\int_{\mathbb R^N}h(y) (u(y))^{\tau}\rho(\frac{y}{R}) dy\right)^{\frac{2}{\tau}}\\
&\leq C\left(R^{-N-\ell}\int_{\mathbb R^N} h(y)(u(y))^{\tau}\rho(\frac{y}{R}) dy\right)^{\frac{2}{\tau}}.
\end{split}
\end{align}
We then split the rest of the proof into two cases.

\medskip
\noindent
\textbf{Case 1}. $q>\frac43.$ Fix $\tau=q.$ By \cite[Lemma 6]{HHM16}, $\tau $ satisfies $(*).$ Using the estimate \eqref{2.6}, we derive from \eqref{3.15} that
\begin{align}\label{3.16}
\int_{B^+_R} U^{\gamma_mk_s}(x,t)t^{1-2s}dxdt&\leq CR^{N+2-2s-\frac{1}{q}\left(\frac{2sq(p+1)}{pq-1}+\frac{ \ell q(p+1)}{pq-1}\right)\gamma_mk_s}.
\end{align}
Suppose now $N<2s+\frac{(2s+\ell)(p+1)}{pq-1}\gamma_0,$ we can choose $k\in [1,k_s]$ so that $\gamma_m$ is sufficiently close to $\gamma_0$ satisfying
$$N-2s-\frac{(2s+\ell)(p+1)}{pq-1}\gamma_m<0.$$ Let $R$ tend to infinity in \eqref{3.16}, we have a contradiction since $u$ is positive. In other words, the system \eqref{1.1} has no stable solution if $N<2s+(2s+\ell)x_0$ where $x_0=\frac{p+1}{pq-1}\gamma_0.$ Moreover, we can adapt the proof of Remark 3 in \cite{Cow13} to show that
\begin{align*}
2t_0^+\frac{p+1}{pq-1}>4,\quad \forall\; p\geq q>1.
  \end{align*}
By \cite[Lemma 6]{HHM16}, $x_0 \geq 2t_0^+\frac{p+1}{pq-1} > 4$. Therefore, if $N\leq 10+4\ell$, \eqref{1.1} has no stable solution
for any $p \geq q> \frac{4}{3}.$

\medskip
\noindent
\textbf{Case 2}. \textbf{$1<q\leq\frac43,$ and $v$ is bounded.} Put $\tau=2.$ By \cite[Lemma 6]{HHM16}, $\tau$ satisfies $(*).$ The following Lemma is crucial to handle this case. It provides an \textit{a priori} integral estimate of $h u^2\rho $ as for $h u^q\rho,$ and then we can proceed as above to achieve the proof of Theorem \ref{main1}.

\begin{lemma}
\label{l.3.2}
Let $(u,v)$ be a stable solution to \eqref{1.1} with $1< q \leq \min(\frac{4}{3}, p)$. Assume that $v$ is bounded and $h$ satisfies $(H)$, there holds
\begin{align}
\label{3.17}
\int_{\mathbb{R}^N} h u^2\rho(\frac{x}{R}) dx \leq CR^{N-\frac{2(p+1)q}{pq-1}-\frac{\ell(q+1)}{pq-1} - \frac{2(2+\ell)(2 - q)}{p+q-2}}, \quad \forall\; R > 0.
\end{align}
\end{lemma}

\noindent{\bf Proof.} Take $\phi$, $\Phi_k$ as above and $\phi_k(x)=\phi \left(\frac{|x|}{k}\right).$ Combining \eqref{3.3} and \eqref{3.4} with $\gamma =2$,  we get
\begin{align*}
\begin{split}
\int_{\mathbb R^{N+1}_+}|\nabla (U\zeta\Phi_k)|^2t^{1-2s}dxdt
=& \kappa_s \int_{\mathbb R^N}h(x) v(x)^pu(x)\rho(x)\phi_k(x)^2dx\\
&+\int_{\mathbb R^{N+1}_+}U^{2}|\nabla (\zeta\Phi_k) |^2 t^{1-2s}dxdt.
\end{split}
\end{align*}
By \eqref{2.7}, there holds
\begin{align*}
\begin{split}
\int_{\mathbb R^{N+1}_+}|\nabla (U\zeta\Phi_k)|^2t^{1-2s}dxdt
\leq& \kappa_s \sqrt{\frac{p+1}{q+1}}\int_{\mathbb R^N}h(x) v(x)^\frac{p-1}{2}u(x)^{\frac{q-1}{2}}u(x)^{2}\rho(x)\phi_k(x)^2dx\\
&+\int_{\mathbb R^{N+1}_+}U^{2}|\nabla (\zeta\Phi_k) |^2 t^{1-2s}dxdt.
\end{split}
\end{align*}
On the other hand, applying \eqref{3.10} with $\gamma_1=2$, we have
\begin{align*}
\begin{split}
 \kappa_s\sqrt{pq}\int_{\mathbb R^N}h(x) v(x)^\frac{p-1}{2}u(x)^{\frac{q-1}{2}}u(x)^{2}\rho(x)\phi_k(x)^2dx
&\leq \int_{\mathbb R^{N+1}_+}|\nabla (U\zeta\Phi_k)|^2t^{1-2s}dx.
\end{split}
\end{align*}
Combining the two last inequalities and using \eqref{3.13} with $\tau=2$,
 \begin{align}\label{Q}
 \begin{split}
 &\left(\sqrt{p q}-\sqrt{\frac{p+1}{q+1}}\right)\kappa_s \int_{\mathbb R^N}h(x) v(x)^\frac{p-1}{2}u(x)^{\frac{q+3}{2}}\rho(x)\phi_k(x)^2dx \leq  \int_{\mathbb{R}^{N}}h(x) u(x)^{2}\rho(x) dx.
\end{split}
\end{align}
As $v$ is bounded, we can use \eqref{2.8} to deduce that  there exists  $C>0$ such that
  \begin{align*}
 \int_{\mathbb R^N}h(x) u(x)^{\frac{p+q+2}{2}}\rho(x)\phi_k(x)^2dx
&\leq C \|v\|_{\infty}^{a} \int_{\mathbb{R}^{N}}h(x) u(x)^{2}\rho(x) dx,
\end{align*}
 where $a:= \frac{(p-q)(p-1)}{2(q+1)}. $ Letting $k\to \infty$ and using Lebesgue's monotone convergence theorem,
 \begin{align}\label{3.18}
 \int_{\mathbb R^N}h(x) u(x)^{\frac{p+q+2}{2}}\rho(x)dx
&\leq C \|v\|_{\infty}^{a} \int_{\mathbb{R}^{N}}h(x) u(x)^{2}\rho(x) dx.
\end{align}
 Denote
 \begin{align*}
J_1 := \int_{\mathbb R^N}h(x) u(x)^{\frac{p+q+2}{2}}\rho(x)dx , \quad  J_2 :=\int_{\mathbb{R}^{N}}h(x) u(x)^{2}\rho(x) dx.
 \end{align*}
As $1 < q \leq \min(p,\,\frac{4}{3}),$ we have $q < 2 < \frac{p+q + 2}{2}$ and a direct calculation yields
 \begin{align*}
2 = q\lambda + \frac{p+q + 2}{2}(1 - \lambda) \quad \mbox{with } \lambda = \frac{ p +q- 2}{p + 2 - q} \in (0, 1).
 \end{align*}
 By H\"older's inequality  and \eqref{3.18}, we get
  \begin{align*}
J_2 \leq J_1^{1 - \l}
\left(\int_{\mathbb{R}^{N}}h(x) u(x)^{q}\rho(x) dx\right)^{\l}
 \leq \left(C\|v\|_{\infty}^{a}J_2\right)^{1 - \l} \left(\int_{\mathbb{R}^{N}}h u^{q}\rho(x) dx\right)^{\l},
 \end{align*}
which implies
\begin{align*}
J_2 = \int_{\mathbb{R}^{N}}h(x) u(x)^{2}\rho(x) dx \leq  C \|v\|_{\infty}^{a\frac{1-\gamma}{\gamma}}\int_{\mathbb{R}^{N}}h u^{q}\rho(x) dx.
 \end{align*}
By scaling argument as above, for $R>0$ we get
\begin{align*}
 &\quad R^{2\left[\frac{2s(p+1)}{pq-1}+\frac{\ell(p+1)}{pq-1}\right]}\int_{\mathbb{R}^{N}}h(Rx) u^{2}(Rx)\rho(x) dx \\&\leq  C \left(R^{\frac{2s(q+1)}{pq-1}+\frac{\ell(q+1)}{pq-1}}\|v(R.)\|_{\infty}\right)^{a\frac{1-\gamma}{\gamma}}\int_{\mathbb{R}^{N}}h(Rx) u^{q}(Rx)\rho(x) dx.
 \end{align*}
Making a change of variables and using \eqref{2.6}, using again the boundedness of $v$, we deduce that
\begin{align*}
 \int_{\mathbb{R}^{N}}h(x) u^{2}(x)\rho(\frac{x}{R}) dx \leq  C R^{N- \frac{2s q(p+1)}{pq-1}-\frac{ \ell(q+1)}{pq-1}+a\frac{1-\gamma}{\gamma}\left[\frac{2s(q+1)}{pq-1}+\frac{\ell(q+1)}{pq-1}\right]-2\left[\frac{2s(p+1)}{pq-1}+\frac{\ell(p+1)}{pq-1}\right]}.
 \end{align*}
A straightforward computation shows that the above exponent of $R$ is the same stated in \eqref{3.17}, so we are done. \qed

\medskip\noindent
\textbf{Proof of Corollary \ref{main2}.} Let $u$ be a stable solution of equation \eqref{1.12}, then $v = u$ verify the system \eqref{1.1} with $p = q$.  Moreover, we have $$t_0^{\pm}= p\pm\sqrt{p^2-p}$$ and
\begin{align*}
  L(\gamma) = \gamma^4 -16p^2\gamma^2 +
32p^2\gamma-16p^2=(\gamma^2+4p(\gamma-1))(\gamma-2t_0^-)(\gamma-2t_0^+).
\end{align*}
As $t_0^+ > p > 1,$ it follows that $2t_0^+$ is the largest root of $L$ as $t_0^+ > p > 1$. Therefore
$$x_0 =\frac{2t_0^+}{p-1}= \frac{2p+2\sqrt{p^2-p}}{p-1}$$ is the largest root of $H.$ Then, applying Theorem \ref{main1}, the result follows immediately. \qed

\bigskip\noindent	
{\bf Acknowledgments}. I would like to thank Professor Dong Ye for many helpful comments.


\begin{thebibliography}{10}


\bibitem{CS07}
{\sc Caffarelli, L., and Silvestre, L.}
\newblock An extension problem related to the fractional {L}aplacian.
\newblock {\em Comm. Partial Differential Equations 32}, 7-9 (2007),
  1245--1260.


\bibitem{CDG14}
{\sc Chen, W., Dupaigne, L., and Ghergu, M.}
\newblock A new critical curve for the Lane-Emden system.
\newblock{\em Discrete Contin. Dyn. Syst. 34}, (2014), 2469--2479.

\bibitem{CW17}
{\sc Chen, W., and Wang, H.}
\newblock Liouville theorems for the weighted Lane-Emden equation with
finite Morse indices.
 \newblock{\em Math. Methods Appl. Sci. 40}, (2017), 4674--4682.

\bibitem{Cow13}
{\sc Cowan, C.}
\newblock Liouville theorems for stable {L}ane-{E}mden systems and biharmonic
  problems.
\newblock {\em Nonlinearity 26}, 8 (2013), 2357--2371.


\bibitem{CF12}
{\sc Cowan, C., and Fazly, M.}
\newblock On stable entire solutions of semi-linear elliptic equations with
  weights.
\newblock {\em Proc. Amer. Math. Soc. 140}, 6 (2012), 2003--2012.

\bibitem{CG14}
{\sc Cowan, C., and Ghoussoub, N.}
\newblock {Regularity of semi-stable solutions to fourth order nonlinear eigenvalue problems on general domains},
\newblock \emph{ Calc. Var. PDE.}, \textbf{49} (2014), 291--305.



 \bibitem{DDG}
{\sc Dancer, E.~N., Du, Y., and Guo, Z.}
\newblock Finite Morse index solutions of an elliptic equation with
supercritical exponent.
 \newblock {\em J. Differ. Equ. 250}, (2011), 3281--3310.

\bibitem{DG13}
{\sc Du, Y., and Guo, Z.}
\newblock Finite Morse-index solutions and asymptotics of weighted nonlinear
elliptic equations.
 \newblock {\em Adv. Differ. Equ. 18}, (2013), 737--768.

\bibitem{DDW17}
{\sc D\'{a}vila, J., Dupaigne, L., and Wei, J.}
\newblock On the fractional {L}ane-{E}mden equation.
\newblock {\em Trans. Amer. Math. Soc. 369}, 9 (2017), 6087--6104.

\bibitem{DMV17}
{\sc Dipierro, S., Medina, M., and Valdinoci, E.}
\newblock Fractional elliptic problems with critical growth in the
whole of $\mathbb{R}^n.$
\newblock {\em N . Appunti. Scuola Normale Superiore di Pisa (Nuova Serie) [Lecture Notes. Scuola Normale
Superiore di Pisa (New Series)], 15. Edizioni della Normale, Pisa}, 2017. viii+152 pp.

\bibitem{Duong17}
{\sc Duong, A.~T.}
\newblock A {L}iouville type theorem for non-linear elliptic systems involving
  advection terms.
\newblock {\em Complex Var. Elliptic Equ. 63}, 12 (2018), 1704--1720.

\bibitem{DuongPham20}
{\sc Duong, A.~T., and Pham, D.~H.}
\newblock Liouville-type Theorem for Fractional Kirchhoff Equations
with Weights.
\newblock {\em Bulletin of the Iranian Mathematical Society},  (2020).

\bibitem{TuanHoang21}
{\sc Duong, A.~T., and Nguyen, V.~H.}
\newblock Liouville type theorems for some fractional elliptic problems
\newblock {\em Nonlinear Anal. 210}, (2021), 112383.
\bibitem{DP17}
{\sc Duong, A.~T., and Phan, Q.~H.}
\newblock Liouville type theorem for nonlinear elliptic system involving
  {G}rushin operator.
\newblock {\em J. Math. Anal. Appl. 454}, 2 (2017), 785--801.



\bibitem{DF10}
{\sc Dupaigne, L., and Farina, A.}
\newblock Stable solutions of {$-\Delta u=f(u)$} in {$\Bbb R^N$}.
\newblock {\em J. Eur. Math. Soc. (JEMS) 12}, 4 (2010), 855--882.

\bibitem{DFS13}
{\sc  Dupaigne, L., Farina, A., and Sirakov, B.}
\newblock {Regularity of the extremal solutions for the Liouville system, in: Geometric Partial Differential Equations},
\newblock \emph{Publications of the Scuola Normale Superiore/CRM Series}, \textbf{15} (2013),  139--144.



\bibitem{Far07}
{\sc Farina, A.}
\newblock On the classification of solutions of the {L}ane-{E}mden equation on
  unbounded domains of {$\R^N$}.
\newblock {\em J. Math. Pures Appl. (9) 87}, 5 (2007), 537--561.




\bibitem{FH19}
{\sc Farina, A., and Hasegawa, S.}
\newblock Liouville-type theorems and existence results for stable solutions to weighted Lane-Emden equations.
\newblock {\em Proceedings of the Royal Society of Edinburgh Section A Mathematics (3) 150}, (2019), 1--13.

\bibitem{FG13}
{\sc Fazly, M., and Ghoussoub, N.}
\newblock On the {H}\'enon-{L}ane-{E}mden conjecture.
\newblock {\em Discrete Contin. Dyn. Syst. 34}, 6 (2014), 2513--2533.

\bibitem{FW16}
{\sc Fazly, M., and Wei, J.}
\newblock On stable solutions of the fractional {H}\'{e}non-{L}ane-{E}mden
  equation.
\newblock {\em Commun. Contemp. Math. 18}, 5 (2016), 1650005, 24.

\bibitem{FW17}
{\sc Fazly, M., and Wei, J.}
\newblock On finite {M}orse index solutions of higher order fractional
  {L}ane-{E}mden equations.
\newblock {\em Amer. J. Math. 139}, 2 (2017), 433--460.

\bibitem{gnw}{\sc Gui, C., Ni, W., and Wang, X.}
\newblock On the stability and instability of positive steady states of a semilinear heat equation in ${\bf R}^n$.
\newblock {\em Comm. Pure Appl. Math. Vol. XLV}, (1992), 1153-1181.

\bibitem{HHM16}
{\sc Hajlaoui, H., Harrabi, A., and Mtiri, F.}
\newblock Liouville theorems for stable solutions of the weighted
  {L}ane-{E}mden system.
\newblock {\em Discrete Contin. Dyn. Syst. 37}, 1 (2017), 265--279.

\bibitem{Har21}
{\sc Harrabi, A.}
\newblock Explicit universal estimate for p-polyharmonic equations via Morse index.
\newblock {\em  arXiv:2105.04058v1} (2021).


\bibitem{Hu15}
{\sc Hu, L.-G.}
\newblock Liouville type results for semi-stable solutions of the weighted
  {L}ane-{E}mden system.
\newblock {\em J. Math. Anal. Appl. 432}, 1 (2015), 429--440.

\bibitem{Hu18}
{\sc Hu, L.-G.}
\newblock Liouville type theorems for stable solutions of the weighted elliptic
  system with the advection term: {$p\ge\vartheta>1$}.
\newblock {\em NoDEA Nonlinear Differential Equations Appl. 25}, 1 (2018), Art.
  7, 30.

\bibitem{HZ16}
{\sc Hu, L.-G., and Zeng, J.}
\newblock Liouville type theorems for stable solutions of the weighted elliptic
  system.
\newblock {\em J. Math. Anal. Appl. 437}, 2 (2016), 882--901.


\bibitem{JL13}
{\sc Jeong, W., and Lee, Y.}
\newblock Stable solutions and finite Morse index solutions of nonlinear elliptic
equations with Hardy potential.
\newblock {\em  Nonlinear Anal. 87}, (2013), 126--145.

\bibitem{jl}
{\sc Joseph, D.D., and Lundgren, T.S.}
\newblock Quasilinear Dirichlet problems driven
by positive sources.
 \newblock\emph{Arch. Rational Mech. Anal. 49},
(1973), 241-269.





\bibitem{Mit96}
{\sc Mitidieri, E.}
\newblock Nonexistence of positive solutions of semilinear elliptic systems in
  {${\mathbb R}^ N$}.
\newblock {\em Differential Integral Equations 9}, 3 (1996), 465--479.

\bibitem{MPO}
{\sc Mitidieri, E., and Pohozaev, S.}
\newblock A priori estimates and the absence of solutions of nonlinear
partial differential equations and inequalities.
\newblock \emph{Tr. Mat. Inst. Steklova 234}, (2001), 1-384.


\bibitem{Mon05}
{\sc Montenegro, M.}
\newblock Minimal solutions for a class of elliptic systems.
\newblock {\em Bull. London Math. Soc. 37}, 3 (2005), 405--416.

\bibitem{MY19}
{\sc Mtiri, F., and Ye, D.}
\newblock Liouville theorems for stable at infinity solutions of {L}ane-{E}mden
  system.
\newblock {\em Nonlinearity 32}, 3 (2019), 910--926.

\bibitem{QS12}
{\sc Quittner, P., and Souplet, P.}
\newblock Symmetry of components for semilinear elliptic systems.
\newblock {\em SIAM J. Math. Anal. 44}, 4 (2012), 2545--2559.

\bibitem{RC19}
{\sc Rahal, B., and Zaidi, C.}
\newblock On the classification of stable solutions of the fractional equation.
\newblock {\em Potential Anal. 50}, 4 (2019), 565--579.


\bibitem{SZ96}
{\sc Serrin, J., and Zou, H.}
\newblock Non-existence of positive solutions of {L}ane-{E}mden systems.
\newblock {\em Differential Integral Equations 9}, 4 (1996), 635--653.


\bibitem{Sou09}
{\sc Souplet, P.}
\newblock The proof of the {L}ane-{E}mden conjecture in four space dimensions.
\newblock {\em Adv. Math. 221}, 5 (2009), 1409--1427.


\bibitem{WY12}
{\sc Wang, C., and Ye, D.}
\newblock Some {L}iouville theorems for {H}\'enon type elliptic equations.
\newblock {\em J. Funct. Anal. 262}, 4 (2012), 1705--1727.

\bibitem{WY13}
 {\sc Wei, J., and Ye, D.}
\newblock Liouville theorems for stable solutions of biharmonic problem,
 \newblock{\em Math. Ann. 356}, (2013), 1599-1612.

\bibitem{YZ19_b}
{\sc Yang, H., and Zou, W.}
\newblock Symmetry of components and {L}iouville-type theorems for semilinear
  elliptic systems involving the fractional {L}aplacian.
\newblock {\em Nonlinear Anal. 180\/} (2019), 208--224.

\end{thebibliography}

\end{document}